\documentclass[a4paper,11pt]{amsart}
\usepackage{amssymb}
\usepackage{graphicx}
\usepackage[all]{xy}

\title[Closed symmetric monoidal structure and flow]{Closed symmetric monoidal structure and flow}
\author[P. Gaucher]{Philippe Gaucher}
\address{Institut de Recherche Math\'ematique Avanc\'ee\\ ULP et
CNRS\\ 7 rue Ren\'e-Descartes\\ 67084 Strasbourg Cedex\\ France}
\email{gaucher@math.u-strasbg.fr}
\urladdr{http://www-irma.u-strasbg.fr/\~{}gaucher/}

\newcommand{\C}{\mathcal{C}}

\newcommand{\de}{\partial}
\newcommand{\p}\times
\renewcommand{\vec}{\overrightarrow}
\renewcommand{\P}{\mathbb{P}}

\newcommand{\be}{\begin{equation}}
\newcommand{\ee}{\end{equation}}
\newcommand{\bea}{\begin{eqnarray}}
\newcommand{\eea}{\end{eqnarray}}
\newcommand{\beas}{\begin{eqnarray*}}
\newcommand{\eeas}{\end{eqnarray*}}

\newtheorem{thm}{Theorem}[section]
\newtheorem{prop}[thm]{Proposition}

\newtheorem{cor}[thm]{Corollary}
\newtheorem{rem}[thm]{Remark}

\newtheorem{defn}[thm]{Definition}

\newtheorem{nota}[thm]{Notation}
\newcommand{\bd}{\begin{defn}}
\newcommand{\ed}{\end{defn}}
\newcommand{\bcd}{\begin{defn}}
\newcommand{\ecd}{\end{defn}}
\newcommand{\bex}{\begin{exmp}}
\newcommand{\eex}{\end{exmp}}
\newcommand{\bp}{\begin{prop}}
\newcommand{\ep}{\end{prop}}
\newcommand{\bth}{\begin{thm}}
\renewcommand{\eth}{\end{thm}}
\newcommand{\br}{\begin{rem}}
\newcommand{\er}{\end{rem}}
\newcommand{\bpf}{\begin{proof}}
\newcommand{\epf}{\end{proof}}

\newcommand{\fl}[1]{\ar@{->}[l]_{#1}}
\newcommand{\fr}[1]{\ar@{->}[r]^-{#1}}
\newcommand{\fd}[1]{\ar@{->}[d]_{#1}}
\newcommand{\fu}[1]{\ar@{->}[u]^{#1}}
\newcommand{\f}[2]{\ar@{->}[#1]|{#2}}
\newcommand{\ff}[2]{\ar@2{->}[#1]|{#2}}
\newcommand{\frr}[1]{\ar@{->}[rr]^{#1}}

\renewcommand{\top}{{\mathbf{Top}}}

\newcommand{\iso}{\cong}
\newcommand{\lp}{\left(}
\newcommand{\rp}{\right)}

\newcommand{\ot}{\otimes}

\newcommand{\vI}{\vec{I}}

\renewcommand{\geq}{\geqslant}

\newcommand{\cattopn}{{\brm{1Cat^{top}_1}}}
\newcommand{\tcattopn}{{\brm{1CAT^{top}_1}}}

\renewcommand{\hom}{[\dtop]}
\newcommand{\homcat}{[\cattopn]}

\def\cartesien{%
  \ar@{-}[]+R+<6pt,-2pt>;[]+RD+<6pt,-6pt>%
  \ar@{-}[]+D+<2pt,-6pt>;[]+RD+<6pt,-6pt>%
}
\def\cocartesien{%
  \ar@{-}[]+L+<-6pt,+2pt>;[]+LU+<-6pt,+6pt>%
  \ar@{-}[]+U+<-2pt,+6pt>;[]+LU+<-6pt,+6pt>%
}

\newcommand{\brm}[1]{\rm{\mathbf{#1}}}

\renewcommand{\top}{{\brm{Top}}}

\newcommand{\dtop}{{\brm{Flow}}}
\newcommand{\dcgtop}{{\brm{Flow}}}

\newcommand{\tdtop}{{\brm{FLOW}}}
\newcommand{\ttop}{{\brm{TOP}}}

\newcommand{\glob}{{\rm{Glob}}}

\newcommand{\liminj}{\varinjlim}
\newcommand{\limproj}{\varprojlim}

\makeatletter

\def\varholim@#1#2{%
  \vtop{\m@th\ialign{##\cr
    \hfil$#1\operator@font holim$\hfil\cr
    \noalign{\nointerlineskip\kern1.5\ex@}#2\cr
    \noalign{\nointerlineskip\kern-\ex@}\cr}}%
}
\def\holimproj{%
  \mathop{\mathpalette\varholim@{\leftarrowfill@\textstyle}}\nmlimits@
}
\def\holiminj{%
  \mathop{\mathpalette\varholim@{\rightarrowfill@\textstyle}}\nmlimits@
}

\newskip\@bigflushglue \@bigflushglue = -100pt plus 1fil

\def\bigcentering{\let\\\@centercr\rightskip\@bigflushglue%
\leftskip\@bigflushglue
\parindent\z@\parfillskip\z@skip}

\makeatother

\begin{document}

\begin{abstract}
The category of flows is not cartesian closed. We construct a closed
symmetric monoidal structure which has moreover a satisfactory
behavior from the computer scientific viewpoint.
\end{abstract}

\maketitle

\tableofcontents

\section{Introduction}

The category of flows was introduced in \cite{model3} as a convenient
framework for the study of higher dimensional automata up to
homotopy. However the category of flows is not cartesian closed
(cf. Proposition~\ref{notcartesian}).  Moreover, the categorical
product of flows is badly behaved from a computer scientific viewpoint
for very simple reasons. Indeed, if $X$ and $Y$ are two flows with non
empty path spaces $\P X$ and $\P Y$, then the path space $\P(X\p Y)$
is isomorphic to $\P X\p \P Y$. In other terms, if $\gamma$ is a
non-constant execution path of $X$ and if $\alpha\in Y^0$ is a state
of $Y$, then $(\gamma,\alpha)$ does not correspond to any non-constant
execution path of $X\p Y$. This problem disappears by considering the
tensor product $X\ot Y$ (cf. Proposition~\ref{deftenseur}) since the
path space $\P(X\ot Y)$ is exactly $(\P X\p Y^0)\sqcup (X^0\p \P Y)
\sqcup (\P X\p \P Y)$ with the composition law that one expects to
find. And it turns out that this new symmetric monoidal structure is
closed.

Section~\ref{reminder} is a short reminder about
flows. Section~\ref{re1cat} recalls the definition of a
\textit{non-contracting topological $1$-category} and 
the construction of the closed monoidal structure $\ot$ made in
\cite{model3}. Section~\ref{tenseurflow} shows that the tensor product
of flows is actually a closed monoidal structure, using
Section~\ref{re1cat}. Section~\ref{explicit} provides a more explicit
way of establishing the closedness of the tensor product of flows. This
new method does not seem to be applicable to the case of non-contracting
topological $1$-categories. Section~\ref{modelmonoidal} proves a
negative result: this new closed monoidal structure together with the
model structure constructed in \cite{model3} does not provide the
category of flows a structure of monoidal model category.

\section{Warning}

This paper belongs to a set of papers corresponding to the
content of ``A Convenient Category for the Homotopy Theory of
Concurrency'' \cite{flow}. Indeed, the latter paper does not seem to
be publishable because of its length. A detailed abstract (in French)
of \cite{flow} can be found in \cite{pgnote1} and \cite{pgnote2}.

\section{Reminder about the category of flows}\label{reminder}

\subsection{Compactly generated topological space}

The category $\top$ of compactly generated topological spaces is
complete, cocomplete and cartesian closed (more details for this
kind of topological spaces in \cite{MR90k:54001,MR2000h:55002},
the appendix of \cite{Ref_wH} and also the preliminaries of
\cite{model3}). Let us denote by $\ttop\left(X,-\right)$ the right
adjoint of the functor $-\p X:\top\longrightarrow \top$. For the
sequel, any topological space will be supposed to be compactly
generated.

Let $n\geq 1$. Let $\mathbf{D}^n$ be the closed $n$-dimensional
disk. Let $\mathbf{S}^{n-1}=\de \mathbf{D}^n$ be the boundary of
$\mathbf{D}^n$ for $n\geq 1$. Notice that $\mathbf{S}^0$ is the
discrete\index{discrete topology} two-point topological space
$\{-1,+1\}$.  Let $\mathbf{D}^0$ be the one-point topological
space. Let $\mathbf{S}^{-1}=\varnothing$ be the empty set.

\subsection{Definition of a flow}

\bd A {\rm flow} $X$ consists of a topological space $\P X$, a
discrete space $X^0$, two continuous maps $s$ and $t$ from $\P X$
to $X^0$ and a continuous and associative map $*:\{(x,y)\in \P
X\p \P X; t(x)=s(y)\}\longrightarrow \P X$ such that
$s(x*y)=s(x)$ and $t(x*y)=t(y)$.  A morphism of flows
$f:X\longrightarrow Y$ consists of a set map
$f^0:X^0\longrightarrow Y^0$ together with a continuous map $\P
f:\P X\longrightarrow \P Y$ such that $f(s(x))=s(f(x))$,
$f(t(x))=t(f(x))$ and $f(x*y)=f(x)*f(y)$. The corresponding
category will be denoted by $\dtop$. \ed

The continuous map $s:\P X\rightarrow X^0$ is called the
\textit{source map}. The continuous map $t:\P X\rightarrow X^0$
is called the \textit{target map}. One can canonically extend
these two maps to the whole underlying topological space
$X^0\sqcup \P X$ of $X$ by setting $s\left(x\right)=x$ and
$t\left(x\right)=x$ for $x\in X^0$.

The topological space $X^0$ is called the
\textit{$0$-skeleton\index{$0$-skeleton}} of $X$. The
$0$-dimensional elements of $X$ are also called \textit{states} or
\textit{constant execution path\index{execution path}s}.

The elements of $\P X$ are called \textit{non constant execution
path\index{execution path}s}.  If $\gamma_1$ and $\gamma_2$ are
two non-constant execution path\index{execution path}s, then
$\gamma_1 *\gamma_2$ is called the
concatenation\index{concatenation} or the composition of
$\gamma_1$ and $\gamma_2$. For $\gamma\in \P X$,
$s\left(\gamma\right)$ is called the \textit{beginning} of
$\gamma$ and $t\left(\gamma\right)$ the \textit{ending} of
$\gamma$.

\begin{nota} For $\alpha,\beta\in X^0$, let $\P_{\alpha,\beta}X$ be the
subspace of $\P X$ equipped the
Kelleyfication\index{Kelleyfication} of the relative topology
consisting of the non-execution path\index{execution path}s of
$X$ with beginning $\alpha$ and with ending $\beta$. \end{nota}

\begin{nota} An element $x\in X^0$ is \textit{achronal} if $x\notin
s\left(\P X\right)\cup t\left(\P X\right)$.  The subspace of
achronal points of $X^0$ is denoted by $Ach\left(X\right)$. The
image of $s$ is denoted by $\P_0 X$ and the image of $t$ by $\P_1
X$. Therefore $Ach\left(X\right)=X^0\backslash\left(\P_0 X \cup
\P_1 X\right)$. \end{nota}

\bd For $X$ a flow, a point $\alpha$ of $X^0$ such that there are
no non-constant execution paths $\gamma$ such that
$t\left(\gamma\right)=\alpha$ (resp.
$s\left(\gamma\right)=\alpha$) is called {\rm initial state}
(resp. {\rm final state}). \ed

\bd\cite{model3}
Let $Z$ be a topological space. Then the {\rm globe} of $Z$ is the
flow $\glob(Z)$ defined as follows: $\glob(Z)^0=\{0,1\}$,
$\P\glob(Z)=Z$, $s=0$, $t=1$ and the composition law is trivial. \ed

\bd\cite{model3}
The {\rm directed segment} $\vI$ is the flow defined as
follows: $\vI^0=\{0,1\}$, $\P\vI=\{[0,1]\}$, $s=0$ and
$t=1$. \ed

\begin{nota} The space $\tdtop(X,Y)$ is the set $\dtop(X,Y)$
equipped with the Kelleyfication of the compact-open topology.
\end{nota}

\bth\cite{model3}\label{lim-colim} The category $\dtop$ is complete and
cocomplete. In particular, a terminal object is the flow $\mathbf{1}$
having the discrete\index{discrete topology} set $\{0,u\}$ as
underlying topological space with $0$-skeleton\index{$0$-skeleton}
$\{0\}$ and path space $\{u\}$. And a initial object is the unique
flow $\varnothing$ having the empty set as underlying topological
space.  \eth

\bd\cite{model3} A morphism of flows $f:X\rightarrow Y$ is said
{\rm synchronized} if and only if it induces a bijection of sets
between the $0$-skeleton\index{$0$-skeleton} of $X$ and the
$0$-skeleton\index{$0$-skeleton} of $Y$. \ed

\section{Reminder about non-contracting topological $1$-categories}\label{re1cat}

\bd \cite{model3}
A {\rm non-contracting topological $1$-category} $X$ is a pair of
compactly generated topological spaces $(X^0,\P X)$ together with
continuous maps $s$, $t$ and $*$ satisfying the same properties as in
the definition of flow except that $X^0$ is not necessarily
discrete.  The corresponding category is denoted by $\cattopn$. \ed

\bth\label{rr}\cite{model3}
The category $\cattopn$ is complete and cocomplete.  The inclusion
functor $\widetilde{\omega}:\dtop\rightarrow
\cattopn$ preserves finite limit\index{limit}s. \eth

\bth\label{ss}\cite{model3} The inclusion functor
$\widetilde{\omega}:\dtop\rightarrow \cattopn$ has a right
adjoint that will be denoted by $\widetilde{D}$. In particular,
this implies that the canonical inclusion functor
$\dtop\rightarrow \cattopn $ preserves colimit\index{colimit}s.
Moreover, one has $\widetilde{D}\circ
\widetilde{\omega}=Id_{\dtop}$ and
\[\limproj_i X_i\iso \limproj_i\widetilde{D}\circ
\widetilde{\omega}\left(X_i\right)\iso \widetilde{D}
\left(\limproj_i \widetilde{\omega}\left(X_i\right)\right).\] \eth

\bp\label{tenseurn}\cite{model3}
Let $X$ and $Y$ be two objects of $\cattopn$. There exists a unique
structure of topological $1$-category
 $X\ot Y$ on the topological space  $X\p Y$ such that
\begin{enumerate}
\item $\left(X\ot Y\right)^0=X^0\p Y^0$ .
\item $\P\left(X\ot Y\right)= \left(\P X\p \P X\right)\sqcup \left(X^0 \p \P Y\right) \sqcup \left(\P X\p Y^0\right)$.
\item $s\left(x,y\right)=\left(sx,sy\right)$, $t\left(x,y\right)=\left(tx,ty\right)$, $\left(x,y\right)*\left(z,t\right)=\left(x*z,y*t\right)$.
\end{enumerate}
\ep

\bth\label{na}\cite{model3} The tensor product\index{tensor product} of
$\cattopn$ is a closed symmetric monoi\-dal structure, that is there
exists a bifunctor \[ \homcat:\cattopn\p \cattopn
\longrightarrow \cattopn\] contravariant with respect to the first
argument and covariant with respect to the second argument such that
one has the natural isomorphism of sets \[\cattopn\left(X\ot
Y,Z\right)\iso \cattopn\left(X,\homcat\left(Y,Z\right)\right)\] for
any topological $1$-categories $X$, $Y$ and $Z$. Moreover, one has the
natural homeomorphism
\[\lp \homcat\left(Y,Z\right)\rp^0\iso \tdtop(Y,Z).\]
\eth

\section{Tensor product of flows}\label{tenseurflow}

\bp\label{notcartesian}\cite{model3}
The category of flows\index{flow}
$\dcgtop$ (as well as the category of
non-contracting\index{non-contracting $1$-category} topological
$1$-categories $\cattopn$) is not cartesian closed. \ep

\bpf We recall here the proof given in \cite{model3}.
If $\dcgtop$ was cartesian closed, then its product would
commute with colimit\index{colimit}. This latter property fails as we
can see with the following example. Let $2_0$ be the flow consisting
of one achronal point $*$. Consider the flows $\vI_v$ and $\vI_w$ such
that $\vI_v\iso \vI$, $\P \vI_v=\{v\}$, $\vI_w\iso \vI$, $\P
\vI_v=\{w\}$. Consider the diagram of $\dtop$
\[\xymatrix{\vI_v & \fl{i_1} {2_0} \fr{i_0} & \vI_w}\]
where $i_0\left(*\right)=0$ and $i_1\left(*\right)=1$. Then the
colimit\index{colimit} $\vI_v.\vI_w$ of this diagram is the flow
representing the concatenation\index{concatenation} of $v$ and
$w$. For any flow $X$, $\P\left(\vI \p X\right)\iso \P X$ as sets and
the paths of $\vI \p X$ are never composable because
$s\left(\P\left(\vI \p X\right)\right)\subset
\{0\}\p X$ and $t\left(\P\left(\vI \p X\right)\right)\subset
\{1\}\p X$. Moreover, the flow $\vI\p2_0$ is the achronal flow
$\{\left(0,*\right),\left(1,*\right)\}$. Therefore the path space
of the colimit\index{colimit} of
\[\xymatrix{{\vI\p\vI_v} & \fl{i_1} {\vI\p2_0} \fr{i_0} & {\vI\p\vI_w}}\]
consists exactly of the two non-composable non-constant execution
path\index{execution path}s $\left([0,1],v\right)$ and
$\left([0,1],w\right)$. On the contrary, the path space of $\vI \p
\left(\vI_v.\vI_w\right)$ consists exactly of  the three
non-composable non-constant execution path\index{execution path}s
$\left([0,1],v\right)$, $\left([0,1],w\right)$ and
$\left([0,1],v*w\right)$. \epf

\bp\label{deftenseur}
Let $X$ and $Y$ be two flows. There exists a unique structure
of flows $X\ot Y$ on the set $X\p Y$ such that
\begin{enumerate}
\item $\left(X\ot Y\right)^0=X^0\p Y^0$
\item $\P \left(X\ot Y\right)= \left(\P X\p \P Y\right)\cup \left(X^0 \p \P Y\right) \cup \left(\P X\p Y^0\right)$
\item $s\left(x,y\right)=\left(sx,sy\right)$, $t\left(x,y\right)=\left(tx,ty\right)$, $\left(x,y\right)*\left(z,t\right)=\left(x*z,y*t\right)$.
\end{enumerate}
Moreover one has $Ach\left(X\ot Y\right)=Ach\left(X\right)\p
Ach\left(Y\right)$ and $\widetilde{\omega}\left(X\ot
Y\right)=\widetilde{\omega}\left(X\right)\ot
\widetilde{\omega}\left(Y\right)$. \ep

\bpf The first part of the statement is clear and is analogue to
Proposition~\ref{tenseurn}. By definition, the following equality
holds: \beas && \P_0\left(X\ot Y\right)=\left(X^0\p \P_0
Y\right)\cup \left(\P_0 X\p Y^0\right)\cup \left(\P_0 X\p \P_0
Y\right)\\ && \P_1\left(X\ot Y\right)=\left(X^0\p \P_1
Y\right)\cup \left(\P_1 X\p Y^0\right)\cup \left(\P_1 X\p \P_1
Y\right) \eeas Therefore \beas \left(X^0\p Y^0\right)\backslash
\P_0\left(X\ot Y\right)&=&\left(X^0\p \left(Y^0\backslash \P_0
Y\right)\right) \cap
\left(\left(X^0\backslash \P_0 X\right)\p Y^0\right) \\&&\cap \left[\left(X^0\p \left(Y^0\backslash \P_0 Y\right)\right)\cup \left(\left(X^0\backslash \P_0 X\right)\p Y^0\right)\right]\\
&=&\left(\left(X^0\backslash \P_0 X\right)\p \left(Y^0\backslash \P_0 Y\right)\right) \\&&\cap \left[\left(X^0\p \left(Y^0\backslash \P_0 Y\right)\right)\cup \left(\left(X^0\backslash \P_0 X\right)\p Y^0\right)\right]\\
&=& \left(X^0\backslash \P_0 X\right)\p \left(Y^0\backslash \P_0
Y\right) \eeas So \beas
Ach\left(X\ot Y\right)&=& \left(X^0\p Y^0\right)\backslash \P_0\left(X\ot Y\right) \cap \left(X^0\p Y^0\right)\backslash \P_1\left(X\ot Y\right)\\
&=& \left(X^0\backslash \P_0 X\right)\p \left(Y^0\backslash \P_0 Y\right)  \cap \left(X^0\backslash \P_1 X\right)\p \left(Y^0\backslash \P_1 Y\right)\\
&=& Ach\left(X\right)\p Ach\left(Y\right) \eeas The equality
$\widetilde{\omega}\left(X\ot
Y\right)=\widetilde{\omega}\left(X\right)\ot
\widetilde{\omega}\left(Y\right)$ comes from
Proposition~\ref{tenseurn}. \epf

\bth The tensor product\index{tensor product} of $\dtop$ is a
closed symmetric monoidal structure, that is there exists a
bifunctor $\hom:\dtop\p \dtop \rightarrow \dtop$ contravariant
with respect to the first argument and covariant with respect to
the second argument such that one has the natural bijection of
sets {\[\dtop\left(X\ot Y,Z\right)\iso
\dtop\left(X,\hom\left(Y,Z\right)\right)\]} for any flows  $X$,
$Y$ and $Z$. \eth

\bpf One has with Theorem~\ref{ss} and Theorem~\ref{na}:
\beas
\dtop\left(X\ot Y,Z\right)&\iso & \cattopn\left(\widetilde{\omega}\left(X\ot Y\right),\widetilde{\omega}\left(Z\right)\right)\\
&\iso &
\cattopn\left(\widetilde{\omega}\left(X\right)\ot \widetilde{\omega}\left(Y\right),\widetilde{\omega}\left(Z\right)\right)\\
&\iso & \cattopn\left(\widetilde{\omega}\left(X\right),\homcat\left(\widetilde{\omega}\left(Y\right),\widetilde{\omega}\left(Z\right)\right)\right)\\
&\iso & \dtop\left(X,\widetilde{D}\circ
\homcat\left(\widetilde{\omega}\left(Y\right),\widetilde{\omega}\left(Z\right)\right)\right)
\eeas so it suffices to set
$\hom\left(Y,Z\right):=\widetilde{D}\circ
\homcat\left(\widetilde{\omega}\left(Y\right),\widetilde{\omega}\left(Z\right)\right)$.
\epf

In particular one has therefore \[\hom\left(\liminj_i
Y_i,Z\right)\iso \limproj_i \hom\left(Y_i,Z\right)\] and
\[\hom\left(Y,\limproj_i Z_i\right)\iso \limproj_i
\hom\left(Y,Z_i\right).\] The first isomorphism is not surprising
because the functor $\widetilde{\omega}$ commutes with
colimit\index{colimit}s by Theorem~\ref{ss}. So one has \beas
\hom\left(\liminj_i Y_i,Z\right)&\iso & \widetilde{D}\circ \homcat\left(\widetilde{\omega}\left(\liminj_i Y_i\right),\widetilde{\omega}\left(Z\right)\right)\\
&\iso & \widetilde{D}\circ \homcat\left(\liminj_i\widetilde{\omega}\left( Y_i\right),\widetilde{\omega}\left(Z\right)\right)\\
&\iso & \widetilde{D}\circ \limproj_i\homcat\left(\widetilde{\omega}\left( Y_i\right),\widetilde{\omega}\left(Z\right)\right)\\
&\iso & \limproj_i\widetilde{D}\circ \homcat\left(\widetilde{\omega}\left( Y_i\right),\widetilde{\omega}\left(Z\right)\right)\\
&\iso & \limproj_i\hom\left(Y_i,Z\right) \eeas A similar
verification does not seem  to be possible for the second
isomorphism because the functor $\widetilde{\omega}$ does not
commute in general with limit\index{limit}s ! In fact something
slightly more complicated happens.

\bp\label{spe} One has the natural isomorphism of flows
\[\widetilde{D}\circ \homcat\left(\widetilde{\omega}\left( Y\right),Z\right)\iso
\widetilde{D}\circ \homcat\left(\widetilde{\omega}\left( Y\right),
\widetilde{\omega}\circ \widetilde{D}\left(Z\right)\right)\]
Moreover one cannot remove the $\widetilde{D}$ on the left from
this isomorphism. \ep

\bpf Indeed \beas
&&\dtop\left(X,\widetilde{D}\circ\homcat\left(\widetilde{\omega}
\left( Y\right),\widetilde{\omega}\circ
\widetilde{D}\left(Z\right)\right)\right)\\
&\iso & \cattopn\left(\widetilde{\omega}\left(X\right),\homcat\left(\widetilde{\omega}\left( Y\right),\widetilde{\omega}\circ \widetilde{D}\left(Z\right)\right)\right)\\
&\iso & \cattopn\left(\widetilde{\omega}\left(X\right)\ot \widetilde{\omega}\left( Y\right),\widetilde{\omega}\circ \widetilde{D}\left(Z\right)\right)\\
&\iso & \cattopn\left(\widetilde{\omega}\left(X\ot Y\right),\widetilde{\omega}\circ \widetilde{D}\left(Z\right)\right)\\
&\iso & \dtop\left(X\ot Y,\widetilde{D}\left(Z\right)\right)\\
&\iso & \cattopn\left(\widetilde{\omega}\left(X\ot Y\right),Z\right)\\
&\iso & \cattopn\left(\widetilde{\omega}\left(X\right)\ot \widetilde{\omega}\left(Y\right),Z\right)\\
&\iso & \cattopn\left(\widetilde{\omega}\left(X\right),\homcat\left(\widetilde{\omega}\left(Y\right),Z\right)\right)\\
&\iso & \dtop\left(X,\widetilde{D}\circ
\homcat\left(\widetilde{\omega}\left(Y\right),Z\right)\right)
\eeas The conclusion follows by Yoneda\index{Yoneda's lemma}'s
lemma. Now suppose that the isomorphism
\[\homcat\left(\widetilde{\omega}\left( Y\right),Z\right)\iso
\homcat\left(\widetilde{\omega}\left( Y\right),
\widetilde{\omega}\circ \widetilde{D}\left(Z\right)\right)\] was
true. Then one would get, by considering the
$0$-skeleton\index{$0$-skeleton}s of the two members,
\[\tcattopn\left(\widetilde{\omega}\left( Y\right),Z\right)\iso
\tcattopn\left(\widetilde{\omega}\left( Y\right),
\widetilde{\omega}\circ \widetilde{D}\left(Z\right)\right)\] for
any object $Z$ of $\cattopn$. So one would have the isomorphisms
of topological spaces
\begin{alignat*}{2}
 &\tdtop\left(Y,\limproj_i Z_i\right)\\
&\iso
\tcattopn\left(\widetilde{\omega}\left(Y\right),\widetilde{\omega}\left(
\limproj_i Z_i\right)\right)\\
&\iso
\tcattopn\left(\widetilde{\omega}\left(Y\right),\widetilde{\omega}\circ
\widetilde{D} \circ \limproj_i \circ
\widetilde{\omega}\left(Z_i\right)\right)&
\hbox{ using Theorem~\ref{rr}}\\
&\iso  \tcattopn\left(\widetilde{\omega}\left(Y\right),\limproj_i \circ \widetilde{\omega}\left(Z_i\right)\right)\\
&\iso  \left(\homcat\left(\widetilde{\omega}\left(Y\right),\limproj_i \circ \widetilde{\omega}\left(Z_i\right)\right)\right)^0\\
&\iso  \limproj_i\left(\homcat\left(\widetilde{\omega}\left(Y\right),\widetilde{\omega}\left(Z_i\right)\right)\right)^0\\
&\iso  \limproj_i \tcattopn\left(\widetilde{\omega}\left(Y\right),\widetilde{\omega}\left(Z_i\right)\right)\\
&\iso  \limproj_i \tdtop\left(Y,Z_i\right) \end{alignat*} which
is known to be false in general because the functor $\tdtop(Y,-)$ does not commute
with all limits\cite{model3}. \epf

Now we can make the verification of the second isomorphism. Indeed
\begin{alignat*}{2} &\hom\left(Y,\limproj_i Z_i\right)\\ &\iso
\widetilde{D}\circ\homcat\left(\widetilde{\omega}\left(Y\right),
\widetilde{\omega}\left(\limproj_i Z_i\right)\right)\\
&\iso
\widetilde{D}\circ\homcat\left(\widetilde{\omega}\left(Y\right),
\widetilde{\omega}\circ\widetilde{D}\circ \limproj_i \circ
\widetilde{\omega}\left(Z_i\right)\right)&
\hbox{ using Theorem~\ref{rr}}\\
&\iso
\widetilde{D}\circ\homcat\left(\widetilde{\omega}\left(Y\right),
\limproj_i\widetilde{\omega}\left(Z_i\right)\right)&
\hbox{ by Proposition~\ref{spe}}\\
&\iso
\widetilde{D}\circ\limproj_i\homcat\left(\widetilde{\omega}\left(Y\right),
\widetilde{\omega}\left(Z_i\right)\right)\\
&\iso
\limproj_i\widetilde{D}\circ\homcat\left(\widetilde{\omega}\left(Y\right),
\widetilde{\omega}\left(Z_i\right)\right)\\
&\iso  \limproj_i \hom\left(Y,Z_i\right).
\end{alignat*}

\section{Explicit construction of the right adjoint}\label{explicit}

This section is devoted to proving the above fact in a
more explicit way.

If $A$ is a topological space and if $Z$ is a flow, let
$\hom\left(\glob\left(A\right),Z\right)$ be the flow defined as
follows:
\begin{itemize}
\item The $0$-skeleton\index{$0$-skeleton} $\hom\left(\glob\left(A\right),Z\right)^0$
is the set \[\dtop\left(\glob\left(A\right),Z\right)\iso
\bigsqcup_{\left(\alpha,\beta\right)\in Z^0\p Z^0}
\top\left(A,\P_{\alpha,\beta}Z\right)\] equipped with the
discrete\index{discrete topology} topology.
\item  The path space  $\P\hom\left(\glob\left(A\right),Z\right)$ is the disjoint sum (in
$\top$)
\[\bigsqcup_{\left(\alpha,\beta,\gamma,\delta\right)\in Z^0\p Z^0\p Z^0\p Z^0}
\left(A,Z\right)_{\alpha,\beta}^{\gamma,\delta}\] where
$\left(A,Z\right)_{\alpha,\beta}^{\gamma,\delta}$ is the pullback
\[
\xymatrix
{\left(A,Z\right)_{\alpha,\beta}^{\gamma,\delta}\cartesien \fr{}
\fd{} & \top\left(A,\P_{\alpha,\beta}Z\right)\p
\P_{\beta,\delta}Z \fd{}
\\ \P_{\alpha,\gamma}Z\p \top\left(A,\P_{\gamma,\delta}Z\right) \fr{} & \ttop\left(A,\P_{\alpha,\delta}Z\right)
}
\]
where $\top\left(-,-\right)$ means the discrete\index{discrete
topology} topology and $\ttop\left(-,-\right)$ the
Kelleyfication\index{Kelleyfication} of the
compact-open\index{compact-open topology} topology.
\item The source map
$\left(A,Z\right)_{\alpha,\beta}^{\gamma,\delta}\rightarrow
\hom\left(\glob\left(A\right),Z\right)^0$ sends an element of
$\left(A,Z\right)_{\alpha,\beta}^{\gamma,\delta}$ to its
projection on $\top\left(A,\P_{\alpha,\beta}Z\right)$ equipped
with the discrete\index{discrete topology} topology.
\item The target map
$\left(A,Z\right)_{\alpha,\beta}^{\gamma,\delta}\rightarrow
\hom\left(\glob\left(A\right),Z\right)^0$ sends an element of
$\left(A,Z\right)_{\alpha,\beta}^{\gamma,\delta}$ to its
projection on $\top\left(A,\P_{\gamma,\delta}Z\right)$ equipped
with the discrete\index{discrete topology} topology.
\item The composition of an element of
$\left(A,Z\right)_{\alpha,\beta}^{\gamma,\delta}$ with an element
of the pullback
\[
\xymatrix
{\left(A,Z\right)_{\gamma,\delta}^{\zeta,\epsilon}\cartesien
\fr{} \fd{} & \top\left(A,\P_{\gamma,\delta}Z\right)\p
\P_{\delta,\epsilon}Z \fd{}
\\ \P_{\gamma,\zeta}Z\p \top\left(A,\P_{\zeta,\epsilon}Z\right) \fr{} & \ttop\left(A,\P_{\gamma,\epsilon}Z\right)
}
\]
is defined as follows. Consider the pullback $P$ of the diagram
\[
\xymatrix
{P \cartesien \fr{}\fd{} & \left(A,Z\right)_{\alpha,\beta}^{\gamma,\delta}\fd{}\\
\left(A,Z\right)^{\zeta,\epsilon}_{\gamma,\delta}\fr{} &
\top\left(A,\P_{\gamma,\delta}Z\right) }
\]
There are canonical continuous maps
\[P\rightarrow \top\left(A,\P_{\alpha,\beta}Z\right)\p \P_{\beta,\delta}Z\p \P_{\delta,\epsilon}Z\]
and
\[P\rightarrow \P_{\alpha,\gamma}Z\p \P_{\gamma,\zeta}Z\p \top\left(A,\P_{\zeta,\epsilon}Z\right)\]
giving rise to the commutative diagram
\[
\xymatrix {P \fr{} \fd{} &
\top\left(A,\P_{\alpha,\beta}Z\right)\p \P_{\beta,\epsilon}Z \fd{}
\\ \P_{\alpha,\gamma}Z\p \top\left(A,\P_{\zeta,\epsilon}Z\right) \fr{} & \ttop\left(A,\P_{\alpha,\epsilon}Z\right)
}
\]
and therefore to a natural continuous map $P\rightarrow
\left(A,Z\right)_{\alpha,\beta}^{\zeta,\epsilon}$. The latter map
yields a natural associative composition law.
\end{itemize}

\bp For any topological space $A$ and $B$, one has the natural
bijection of sets {\[\dtop\left(\glob\left(B\right)\ot
\glob\left(A\right),Z\right)\iso
\dtop\left(\glob\left(B\right),\hom\left(\glob\left(A\right),Z\right)\right).\]}
\ep

\bpf First of all, one has to calculate the tensor
product\index{tensor product} $\glob\left(B\right)\ot
\glob\left(A\right)$. The latter flow can be conventionally
represented as follows:
\[
\xymatrix{ \left(1,0\right)\ar@{->}[rr]^-{\{1\}\p A}&& \left(1,1\right) \\&&\\
\left(0,0\right)\ar@{->}[uu]^-{B\p\{0\}} \ar@{->}[rruu]^-{B\p A}
\ar@{->}[rr]^-{\{0\}\p A} &&
\left(0,1\right)\ar@{->}[uu]_-{B\p\{1\}} }
\]
where the vertices are the elements of the
$0$-skeleton\index{$0$-skeleton} of $\glob\left(B\right)\ot
\glob\left(A\right)$ and where the labels of the arrows between
them are the components of the path space of
$\glob\left(B\right)\ot \glob\left(A\right)$. Let
$\left(B,A,Z\right)_{\alpha,\beta}^{\gamma,\delta}$ be the
elements $f$ of \[\dtop\left(\glob\left(B\right)\ot
\glob\left(A\right),Z\right)\] such that
$f\left(0,0\right)=\alpha$, $f\left(0,1\right)=\beta$,
$f\left(1,0\right)=\gamma$ and $f\left(1,1\right)=\delta$. It is
helpful to notice that the locations of $\alpha$, $\beta$,
$\gamma$ and $\delta$ in the expression
$\left(B,A,Z\right)_{\alpha,\beta}^{\gamma,\delta}$ correspond to
the locations of $\left(0,0\right)$, $\left(0,1\right)$,
$\left(1,0\right)$ and $\left(1,1\right)$ in the above diagram.
Then by definition of a morphism of flows, one has the pullback
{\scriptsize
\[
\xymatrix
{\left(B,A,Z\right)_{\alpha,\beta}^{\gamma,\delta}\cartesien
\fr{} \fd{} & \top\left(\{0\}\p A,\P_{\alpha,\beta}Z\right)\p
\top\left(B\p\{1\},\P_{\beta,\delta}Z\right) \fd{}
\\ \top\left(B\p\{0\},\P_{\alpha,\gamma}Z\right)\p \top\left(\{1\}\p A,\P_{\gamma,\delta}Z\right) \fr{} & \top\left(B\p A,\P_{\alpha,\delta}Z\right)
}
\] }
so one has the pullback
\[
\xymatrix
{\left(B,A,Z\right)_{\alpha,\beta}^{\gamma,\delta}\cartesien
\fr{} \fd{} & \top\left(A,\P_{\alpha,\beta}Z\right)\p
\top\left(B,\P_{\beta,\delta}Z\right) \fd{}
\\ \top\left(B,\P_{\alpha,\gamma}Z\right)\p \top\left(A,\P_{\gamma,\delta}Z\right) \fr{} & \top\left(B,\ttop\left(A,\P_{\alpha,\delta}Z\right)\right)
}
\]
In the other hand, the set
$\dtop\left(\glob\left(B\right),\hom\left(\glob\left(A\right),Z\right)\right)$
is in natural bijection with
\[\bigsqcup_{u,v\in\dtop\left(\glob\left(A\right),Z\right)}
\top\left(B,\P_{u,v}\hom\left(\glob\left(A\right),Z\right)\right)\]
A pair $\left(u,v\right)\in
\dtop\left(\glob\left(A\right),Z\right)\p
\dtop\left(\glob\left(A\right),Z\right)$ is determined by
$u\left(0\right)=\alpha$, $u(1)=\beta$, $v\left(0\right)=\gamma$,
$v(1)=\delta$ and by the continuous map $\P u\in
\top\left(A,\P_{\alpha,\beta}Z\right)$ and $\P v\in
\top\left(A,\P_{\gamma,\delta}Z\right)$. Then one has a natural
bijection of sets between
\[\dtop\left(\glob\left(B\right),\hom\left(\glob\left(A\right),Z\right)\right)\] and the disjoint sum over
\[\left(\alpha,\beta,\gamma,\delta,m,n\right)\in Z^0\p Z^0\p Z^0\p Z^0\p \top\left(A,\P_{\alpha,\beta}Z\right)\p
\top\left(A,\P_{\gamma,\delta}Z\right)\] of elements $f\in
\top\left(B,\left(A,Z\right)_{\alpha,\beta}^{\gamma,\delta}\right)$
such that such that the composite of $f$ with the canonical
projection map
$\left(A,Z\right)_{\alpha,\beta}^{\gamma,\delta}\rightarrow
\top\left(A,\P_{\alpha,\beta}Z\right)$ is the constant map $m$
and such that the composite of $f$ with the canonical projection
map $\left(A,Z\right)_{\alpha,\beta}^{\gamma,\delta}\rightarrow
\top\left(A,\P_{\gamma,\delta}Z\right)$ is the constant map $n$
hence the result. \epf

We need to recall the following theorems for the sequel:

\bth \label{point-globe}\cite{model3} Any flow is the colimit\index{colimit}
in $\dtop$ of points and globe\index{globe}s in a canonical way,
i.e. there exists for any flow $X$ a diagram
$\mathbb{D}\left(X\right)$ of flows containing only points,
globe\index{globe}s and concatenation\index{concatenation}s of
globe\index{globe}s such that the mapping $X\mapsto
\mathbb{D}\left(X\right)$ is functorial and such that $X\iso
\liminj \mathbb{D}\left(X\right)$ in a canonical way. \eth

\begin{cor} \label{reduction}\cite{model3}
Let $P\left(X\right)$ be a statement depending on a flow $X$ and
satisfying the following property: if
$D:\mathcal{I}\rightarrow \dtop$ is a diagram of flows such that
for any object $i$ of $\mathcal{I}$,
$P\left(D\left(i\right)\right)$ holds, then $P\left(\liminj
D\right)$ holds. Then the following assertions are
equivalent:
\begin{enumerate}
\item[(i)] The statement $P\left(X\right)$ holds for any flow $X$ of $\dtop$.
\item[(ii)] The statements $P\left(\{*\}\right)$ and $P\left(\glob\left(Z\right)\right)$ hold for any
object $Z$ of $\top$.
\end{enumerate}
\end{cor}

\bth\label{hinterne} Let $Y$ be a flow. Then the functor $-\ot
Y:\dtop\rightarrow \dtop$ has a right adjoint denoted by
$\hom\left(Y,-\right)$. So the tensor product\index{tensor
product} of flows is a closed symmetric monoidal structure. In
particular one has the natural isomorphisms of flows
\[\hom\left(\liminj_i Y_i,Z\right)\iso \limproj_i \hom\left(Y_i,Z\right)\]
and
\[\hom\left(Y,\limproj_i Z\right)\iso \hom\left(Y_i,Z_i\right).\]
\eth

\bpf Let $\hom\left(\{*\},Z\right):=Z$ so that for any flow $X$
one has \[\dtop\left(X\ot \{*\},Z\right)\iso
\dtop\left(X,Z\right)\] and let
$\hom\left(\glob\left(A\right),Z\right)$ be defined as above for
any flow $Z$. Using Theorem~\ref{point-globe}, let $Y=\liminj
\mathbb{D}\left(Y\right)$. Then let
\[\hom\left(Y,Z\right):= \limproj_i \hom\left(\mathbb{D}\left(Y\right)(i),Z\right)\]
Then $\hom\left(Y,Z\right)^0\iso \limproj_i
\hom\left(\mathbb{D}\left(Y\right)(i),Z\right)^0$ as set
therefore the following natural bijections of sets hold
\[\hom\left(Y,Z\right)^0\iso \limproj_i \dtop\left(\mathbb{D}\left(Y\right)(i),Z\right)\iso \dtop\left(Y,Z\right).\]
Then for any topological space $B$, one has \beas
\dtop\left(\glob\left(B\right)\ot Y,Z\right)&\iso & \dtop\left(\glob\left(B\right)\ot \left(\liminj \mathbb{D}\left(Y\right)\right),Z\right)\\
&\iso & \limproj \dtop\left(\glob\left(B\right)\ot \mathbb{D}\left(Y\right),Z\right)\\
&\iso & \limproj \dtop\left(\glob\left(B\right),\hom\left(\mathbb{D}\left(Y\right),Z\right)\right)\\
&\iso & \dtop\left(\glob\left(B\right),\limproj\hom\left(\mathbb{D}\left(Y\right),Z\right)\right)\\
&\iso & \dtop\left(\glob\left(B\right),\hom\left(Y,Z\right)\right)
\eeas Moreover one has {
\[\dtop\left(\{*\}\ot Y,Z\right)\iso \dtop\left(Y,Z\right)\iso \hom\left(Y,Z\right)^0\iso \dtop\left(\{*\},\hom\left(Y,Z\right)\right)\]}
so the natural isomorphism $\dtop\left(X\ot Y,Z\right)\iso
\dtop\left(X,\hom\left(Y,Z\right)\right)$ holds if $X$ is a point
or a globe\index{globe}. We could be tempted to concluding that
the proof is complete using Corollary~\ref{reduction}. However
this would not be correct because we do not know yet that
$\limproj_i \hom\left(Y_i,Z\right)\iso \hom\left(\liminj
Y_i,Z\right)$ for any diagram of flows $i\mapsto Y_i$ ! Let $B_1$
and $B_2$ be two topological spaces. Then the set
$\dtop\left(\left(\glob\left(B_1\right)*\glob\left(B_2\right)\right)\ot
Y,Z\right)$ is the pullback {
\[
\xymatrix{
\dtop\left(\left(\glob\left(B_1\right)*\glob\left(B_2\right)\right)\ot Y,Z\right)\cartesien \fr{}\fd{} & \dtop\left(\glob\left(B_2\right)\ot Y,Z\right)\fd{}\\
\dtop\left(\glob\left(B_1\right)\ot Y,Z\right) \fr{} &
\dtop\left(\{*\}\ot Y,Z\right)}
\]}
since the tensor product\index{tensor product} of flows commutes
with colimit\index{colimit}s. So the set
\[\dtop\left(\left(\glob\left(B_1\right)*\glob\left(B_2\right)\right)\ot Y,Z\right)\]
is the pullback {
\[
\xymatrix{
\dtop\left(\left(\glob\left(B_1\right)*\glob\left(B_2\right)\right)\ot Y,Z\right)\cartesien \fr{}\fd{} & \dtop\left(\glob\left(B_2\right),\hom\left(Y,Z\right)\right)\fd{}\\
\dtop\left(\glob\left(B_1\right),\hom\left(Y,Z\right)\right)
\fr{} & \dtop\left(\{*\},\hom\left(Y,Z\right)\right)}
\]}
by the preceding calculations of this proof. So the natural
bijection of sets
{\[\dtop\left(\left(\glob\left(B_1\right)*\glob\left(B_2\right)\right)\ot
Y,Z\right)\iso
\dtop\left(\left(\glob\left(B_1\right)*\glob\left(B_2\right)\right),\hom\left(Y,Z\right)\right)\]}
holds. Now we can conclude using Theorem~\ref{point-globe} by
\beas
\dtop\left(X\ot Y,Z\right)&\iso & \dtop\left(\left(\liminj \mathbb{D}\left(X\right)\right)\ot Y,Z\right)\\
&\iso & \limproj \dtop\left(\mathbb{D}\left(X\right)\ot Y,Z\right)\\
&\iso & \limproj \dtop\left(\mathbb{D}\left(X\right),\hom\left(Y,Z\right)\right)\\
&\iso & \dtop\left(\liminj \mathbb{D}\left(X\right),\hom\left(Y,Z\right)\right)\\
&\iso & \dtop\left(X,\hom\left(Y,Z\right)\right) \eeas \epf

Notice that in general, the topological spaces \[\tdtop\left(X\ot
Y,Z\right)\] and \[\tdtop(X,\hom(Y,Z))\]  are not homeomorphic !
Indeed for $X=\{*\}$, then
\[\tdtop\left(X\ot Y,Z\right)\iso \tdtop\left(Y,Z\right)\] and
\[\tdtop\left(X,\hom\left(Y,Z\right)\right)\iso\hom\left(Y,Z\right)^0\iso
\dtop\left(Y,Z\right)\] which is always discrete\index{discrete
topology}. However one has

\bp Let $B$ be a topological space. Then for any flow $Y$ and
$Z$, there is a canonical homeomorphism
\[\tdtop\left(\glob\left(B\right)\ot Y,Z\right)\iso \tdtop\left(\glob\left(B\right),\hom\left(Y,Z\right)\right).\]
\ep

\bpf If $Y=\{*\}$, then {\[\tdtop\left(\glob\left(B\right)\ot
Y,Z\right)\iso\tdtop\left(\glob\left(B\right),\hom\left(Y,Z\right)\right)\]}

It then suffices to show the homeomorphism {
\[\tdtop\left(\glob\left(B\right)\ot \glob\left(A\right),Z\right)\iso\tdtop\left(\glob\left(B\right),\hom\left(\glob\left(A\right),Z\right)\right).\]}
to be able to conclude with Corollary~\ref{reduction} and
Theorem~\ref{hinterne}.

There is an inclusion of topological spaces
{\[\tdtop\left(\glob\left(B\right)\ot
\glob\left(A\right),Z\right)\rightarrow
\ttop\left(\glob\left(B\right)\ot \glob\left(A\right),Z\right)\]}
so if $\left(B,A,Z\right)_{\alpha,\beta}^{\gamma,\delta}$ is
equipped with the Kelleyfication\index{Kelleyfication} of the
relative topology induced by \[\tdtop\left(\glob\left(B\right)\ot
\glob\left(A\right),Z\right)\] one has the pullback in $\top$:
{\scriptsize
\[
\xymatrix
{\left(B,A,Z\right)_{\alpha,\beta}^{\gamma,\delta}\cartesien
\fr{} \fd{} & \ttop\left(\{0\}\p A,\P_{\alpha,\beta}Z\right)\p
\ttop\left(B\p\{1\},\P_{\beta,\delta}Z\right) \fd{}
\\ \ttop\left(B\p\{0\},\P_{\alpha,\gamma}Z\right)\p \ttop\left(\{1\}\p A,\P_{\gamma,\delta}Z\right) \fr{} & \ttop\left(B\p A,\P_{\alpha,\delta}Z\right)
}
\]}
In the other hand, the topological space
\[\tdtop\left(\glob\left(B\right),\hom\left(\glob\left(A\right),Z\right)\right)\] is
isomorphic to
\[\bigsqcup_{u,v\in\dtop\left(\glob\left(A\right),Z\right)}
\ttop\left(B,\P_{u,v}\hom\left(\glob\left(A\right),Z\right)\right)\]
As above,  the pairs $\left(u,v\right)\in
\dtop\left(\glob\left(A\right),Z\right)\p
\dtop\left(\glob\left(A\right),Z\right)$ are determined by
$u\left(0\right)=\alpha$, $u(1)=\beta$, $v\left(0\right)=\gamma$,
$v(1)=\delta$ and by the continuous map $\P u\in
\top\left(A,\P_{\alpha,\beta}Z\right)$ and $\P v\in
\top\left(A,\P_{\gamma,\delta}Z\right)$. Then one has a natural
isomorphism of topological spaces  between
\[\tdtop\left(\glob\left(B\right),\hom\left(\glob\left(A\right),Z\right)\right)\] and the disjoint sum over
\[\left(\alpha,\beta,\gamma,\delta,m,n\right)\in Z^0\p Z^0\p Z^0\p Z^0\p \top\left(A,\P_{\alpha,\beta}Z\right)\p
\top\left(A,\P_{\gamma,\delta}Z\right)\] of elements $f\in
\ttop\left(B,\left(A,Z\right)_{\alpha,\beta}^{\gamma,\delta}\right)$
such that such that the composite of $f$ with the canonical
projection map
$\left(A,Z\right)_{\alpha,\beta}^{\gamma,\delta}\rightarrow
\ttop\left(A,\P_{\alpha,\beta}Z\right)$ is the constant map $m$
and such that the composite of $f$ with the canonical projection
map $\left(A,Z\right)_{\alpha,\beta}^{\gamma,\delta}\rightarrow
\ttop\left(A,\P_{\gamma,\delta}Z\right)$ is the constant map $n$
hence the result. \epf

\section{Model structure of $\dtop$ and tensor product of flows}\label{modelmonoidal}

Some useful references for the notion of model\index{model
category} category are \cite{MR99h:55031,MR2001d:55012}. See also
\cite{ref_model1,ref_model2}.

\bth \cite{model3}\label{model}  The category of flows
can be given a model structure such that:
\begin{enumerate}
\item The weak equivalences  are the weak S-homotopy equivalences.
\item The fibrations are the continuous maps satisfying the RLP
(right lifting property) with respect
to the morphisms $\glob(\mathbf{D}^n)\longrightarrow \glob([0,1]\p \mathbf{D}^n)$ for $n\geq
0$. The fibration are exactly the morphisms of flows $f:X\longrightarrow
Y$ such that $\P f:\P X\longrightarrow \P Y$ is a Serre fibration of $\top$.
\item The cofibrations are the morphisms satisfying the LLP (left lifting property) with
respect to any map satisfying the RLP with respect to the morphisms
$\glob(\mathbf{S}^{n-1})\longrightarrow \glob(\mathbf{D}^n)$ with $\mathbf{S}^{-1}=\varnothing$ and for
$n\geq 0$ and with respect to the morphisms $\varnothing\longrightarrow
\{0\}$ and $\{0,1\}\longrightarrow \{0\}$.
\item Any flow is fibrant.
\end{enumerate}
\eth

\bp\label{ortho}\cite{model3} Let $U$ and $V$ be two topological spaces. Let
$\phi:U\rightarrow V$ be a continuous map. Let $f:X\rightarrow Y$
be a morphism of flows and let $\P f:\P X\rightarrow \P Y$ the
corresponding continuous maps between the two path spaces. Assume
that at least one of the following conditions holds:
\begin{enumerate}
\item both topological spaces $U$ and $V$ are connected
\item the topological space $V$ is connected and the morphism of flows
$f$ is synchronized.
\end{enumerate}
Then the following conditions are equivalent:
\begin{enumerate}
\item the morphism $f$ satisfies the RLP\index{right lifting property} with respect to the
morphism of flows $\glob(\phi):\glob(U)\rightarrow \glob(V)$
\item the continuous map $\P f$ satisfies the RLP\index{right lifting property} with respect to
the continuous map $\phi:U\rightarrow V$.
\end{enumerate}
\ep

The following definition is an adaptation of the more general
notion of monoidal model\index{model category} category (cf.
\cite{MR99h:55031}).

If $(\C,\otimes)$ is a monoidal category, if $f:U\rightarrow V$
and $g:W\rightarrow X$ are two morphisms of $\C$, then let
\[f\square g:P(f,g)=(V\otimes W)\sqcup_{U\otimes W} (U\otimes X)\rightarrow V\otimes X\]

\bd Let $\C$ be a cofibrantly generated model\index{model
category} category equipped with a closed symmetric monoidal
structure $\otimes$. Let $I$ be the set of generating
cofibrations and let $J$ be the set of generating acyclic
cofibrations. Then $\C$ together with $\otimes$ is a
\textit{monoidal model\index{model category} category} if the
following conditions hold:
\begin{enumerate}
\item The monoidal structure $\otimes$ is a Quillen bifunctor, which means
here that any morphism of $I\square I$ is a cofibration and that
any morphism of $I\square J$ and $J\square I$ is an acyclic
cofibration.
\item Let $q:Q S\rightarrow S$ be the cofibrant replacement for the unit
$S$ of $\otimes$. Then the natural morphism $q\otimes Id_X:Q
S\otimes X\rightarrow S\otimes X$  is a weak equivalence for any
$X$.
\end{enumerate}
\ed

Unfortunately, one has:

\bp The model category $\dtop$ together with the closed
mon\-oi\-dal structure $\otimes$ does not satisfy the axioms of
monoidal model category. \ep

\bpf Consider the two cofibrations of flows $f:\{0,1\}\rightarrow
\{0\}$ and $g:\glob(\mathbf{S}^{n-1})\rightarrow
\glob(\mathbf{D}^n)$ for $n\geq 1$. Then {\small\beas
&&P(f,g)\\&&=\{0\}\otimes \glob(\mathbf{S}^{n-1})
\sqcup_{\{0,1\}\otimes \glob(\mathbf{S}^{n-1})}
\{0,1\}\otimes \glob(\mathbf{D}^n)\\
&&\iso \{0\}\otimes \glob(\mathbf{S}^{n-1}) \sqcup_{\{0\}\otimes
\glob(\mathbf{S}^{n-1})\sqcup\{1\}\otimes \glob(\mathbf{S}^{n-1})}
\left(\{0\}\otimes \glob(\mathbf{D}^n)\sqcup\{1\}\otimes \glob(\mathbf{D}^n)\right)\\
&&\iso \{0\}\otimes \glob(\mathbf{D}^n)\sqcup_{\glob(\mathbf{S}^{n-1})}\{1\}\otimes \glob(\mathbf{D}^n)\\
&&\iso \glob(\mathbf{D}^n \sqcup_{\mathbf{S}^{n-1}} \mathbf{D}^n)\\
&& \iso \glob(\mathbf{S}^n) \eeas} If
$\glob(\mathbf{S}^n)\rightarrow \glob(\mathbf{D}^n)$ was a
cofibration for the model structure of $\dtop$, then by
Proposition~\ref{ortho}, and since $\mathbf{S}^1$ and
$\mathbf{D}^1$ are connected, then $\mathbf{S}^1\rightarrow
\mathbf{D}^1$ would be a cofibration for the model structure of
$\top$. Contradiction. \epf


\begin{thebibliography}{Gau03b}

\bibitem[Bro88]{MR90k:54001}
R.~Brown.
\newblock {\em Topology}.
\newblock Ellis Horwood Ltd., Chichester, second edition, 1988.
\newblock A geometric account of general topology, homotopy types and the
  fundamental groupoid.

\bibitem[DHK97]{ref_model1}
B.~Dwyer, P.~S. Hirschhorn, and D.~Kan.
\newblock Model categories and more general abstract homotopy theory.
\newblock available at http://www-math.mit.edu/\~{}psh/, March 1997.

\bibitem[Gau]{flow}
P.~Gaucher.
\newblock {A Convenient Category for The Homotopy Theory of Concurrency}.
\newblock arXiv:math.AT/0201252.

\bibitem[Gau03a]{model3}
P.~Gaucher.
\newblock {A model category for the homotopy theory of concurrency}, 2003.
\newblock arXiv:math.AT/0308054.

\bibitem[Gau03b]{pgnote1}
P.~Gaucher.
\newblock Concurrent process up to homotopy ({I}).
\newblock {\em C. R. Acad. Sci. Paris Ser. I Math.}, 336(7):593--596, 2003.
\newblock French.

\bibitem[Gau03c]{pgnote2}
P.~Gaucher.
\newblock Concurrent process up to homotopy ({II}).
\newblock {\em C. R. Acad. Sci. Paris Ser. I Math.}, 336(8):647--650, 2003.
\newblock French.

\bibitem[GJ99]{MR2001d:55012}
P.~G. Goerss and J.~F. Jardine.
\newblock {\em Simplicial homotopy theory}.
\newblock Birkh\"auser Verlag, Basel, 1999.

\bibitem[Hir01]{ref_model2}
P.~S. Hirschhorn.
\newblock Localization of model categories.
\newblock available at http://www-math.mit.edu/\~{}psh/, October 2001.

\bibitem[Hov99]{MR99h:55031}
M.~Hovey.
\newblock {\em Model categories}.
\newblock American Mathematical Society, Providence, RI, 1999.

\bibitem[Lew78]{Ref_wH}
L.~G. Lewis.
\newblock {\em The stable category and generalized Thom spectra}.
\newblock PhD thesis, University of Chicago, 1978.

\bibitem[May99]{MR2000h:55002}
J.~P. May.
\newblock {\em A concise course in algebraic topology}.
\newblock University of Chicago Press, Chicago, IL, 1999.

\end{thebibliography}
\end{document}